\input epsf
\epsfverbosetrue

\input amssym.def
\line{}\voffset 0truein
\hoffset 0truein

\baselineskip=14pt
\def\bs{\bigskip}
\def\ms{\medskip}
\def\ss{\smallskip}

\def\ni{\noindent}

\def\<{\langle}
\def\>{\rangle}
\def\a{\alpha}
\def\b{\beta}
\def\g{\gamma}

\def \={ {\buildrel \cdot \over =}}
\def\qed{{$\vrule height4pt depth0pt width4pt$}}

\centerline{\bf ALEXANDER GROUPS OF LONG VIRTUAL KNOTS}\ss

\centerline{Daniel S. Silver and Susan G. Williams} \bs

\ms
\footnote{} {Both authors partially supported by NSF grant
DMS-0071004.}
\footnote{}{2000 {\it Mathematics Subject Classification.}  
Primary 57M25; Secondary 20F05, 20F34.}

{\narrower {\bf ABSTRACT:} Alexander group systems for virtual long knots are defined and used to show that any virtual knot is the closure of infinitely many long virtual knots. Manturov's result that there exists a pair of long virtual knots that do not commute is reproved. \bs


\ni {\bf 1. Introduction.} From a combinatorial perspective, a (classical) knot is an equivalence class of generic immersed circles in the plane with over/under information at each double point. Two such diagrams are equivalent if one can be converted into the other by a finite sequence of Reidemeister moves. 

In 1997 Kauffman introduced virtual knot diagrams by allowing a new sort of crossing, called a {\it virtual crossing}, indicated by a small circle surrounding the site. After augmenting the set of Reidemeister moves, allowing certain deformations  involving virtual crossings, Kauffman defined a {\it virtual knot} to be an equivalence class of virtual diagrams. Details may be found in [K97], [K99]. 

Two classical knots are equivalent under generalized Reidemeister moves if and only they are equivalent under the classical moves [GPV00]. Consequently, virtual knot theory is an extension of the classical theory. Many invariants of classical knots (e.g., the knot group $\pi_\ell$) extend.

A  {\it long virtual knot diagram} is a generic immersion $f: {\Bbb R} \to {\Bbb R}^2$ with over/under/\break virtual crossing information at each double point, and such that $f(x) = (x,0)$ for $|x|$ sufficiently large. Two diagrams are equivalent if one can be transformed into the other by a finite sequence of generalized Reidemeister moves. As before, a {\it long virtual knot} is an equivalence class of diagrams.  In \S2 we extend the Alexander group of a virtual knot, introduced in [SW01], to long virtual knots.

Long virtual knots form a semigroup under concatenation. An example of two long virtual knots that do not commute was given by Manturov [M04] using quandle modules and representations. We give an alternative proof of noncommutativity in \S3 using Alexander groups. 

There is a natural map from the semigroup of long virtual knots to the set of virtual knots given by joining the left and right ends of long virtual knot diagrams. We denote this ``closing map" by $\kappa$. As noted in [GPV00], restricted to the classical category, $\kappa$ is a bijection. The general situation is more complicated. We prove in \S4 that the preimage $\kappa^{-1}(k)$ of any virtual knot $k$ consists of infinitely many distinct long virtual knots. \bs


\ni {\bf 2. Alexander groups.} Let $k$ be an oriented virtual knot with diagram $D$.  The {\it extended Alexander group} $\tilde {\cal A}_k$ was defined in [SW01] (see also [SW03]). We review its definition and key properties for the convenience of the reader. 

In order to define the extended Alexander group, we regard each overcrossing arc of $D$ as a union of two arcs joined at the point of overcrossing. We disregard virtual crossings when determining arcs. The group $\tilde {\cal A}_k$ is described by families of generators $a_{\bf n}, b_{\bf n}, c_{\bf n}\dots$, each indexed by vectors ${\bf n}\in {\Bbb Z}^2$, corresponding to arcs. To each classical crossing we associate an indexed family of relations of the form $a_{\bf n}b_{{\bf n}+u}= c_{\bf n}a_{{\bf n}+u},\ a_{{\bf n}+v} = d_{\bf n}$, where $n\in {\Bbb Z}^2$ (see Figure 1). Here $u, v$ are the standard basis elements of
${\Bbb Z}^2$.  For notational convenience we denote generator families by $a,b,c,\ldots$, and we denote a family of relations  $a_{\bf n}b_{{\bf n}+u}= c_{\bf n}a_{{\bf n}+u}$ by $ab^{u} = c a^{u}$. Similarly, the family $a_{{\bf n}+v} = d_{\bf n}$ is written $a^v = d$. In this notation, exponents $(j,k)\in {\Bbb Z}^2$ are written multiplicatively as $u^jv^k$. \bs

\epsfxsize=2truein
\centerline{\epsfbox{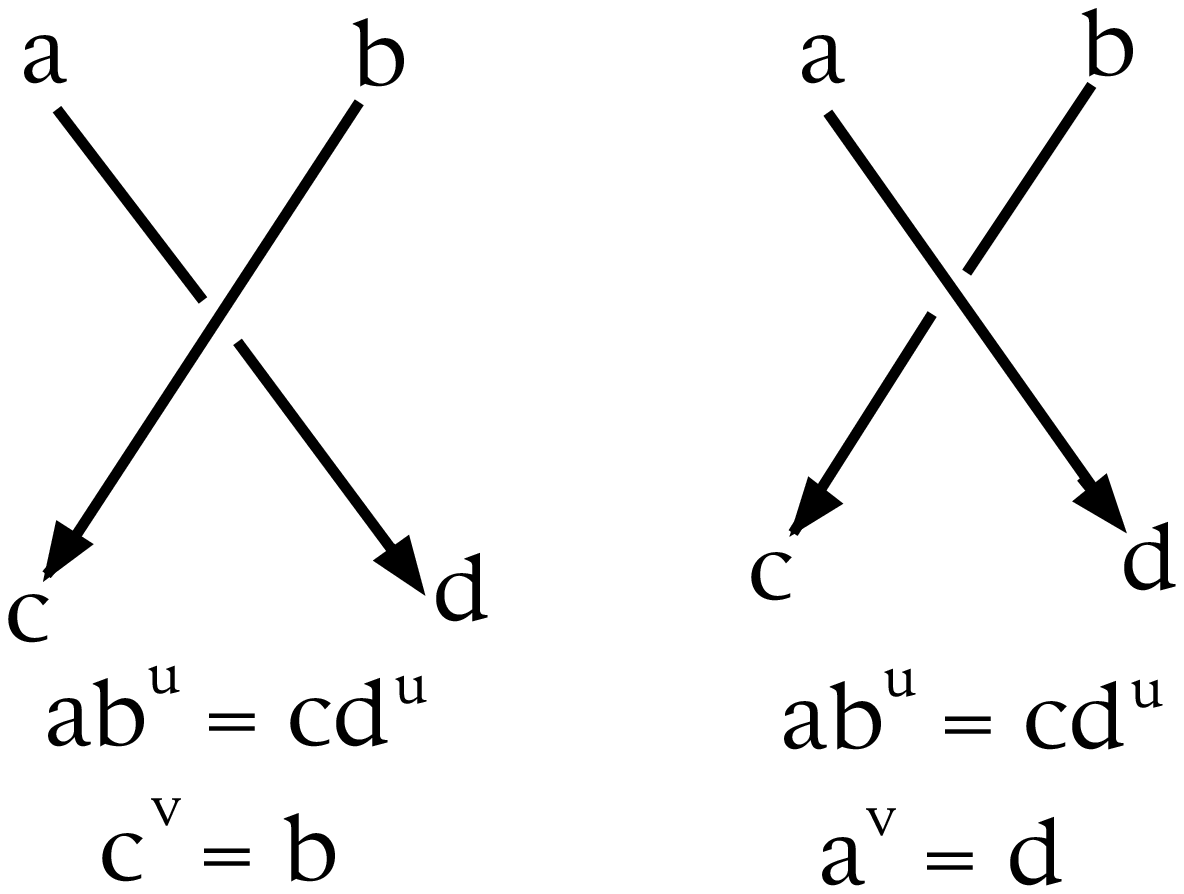}} \bs

\centerline{{\bf Figure 1:} Relation for extended Alexander group}\bs

The extended Alexander group is an example of a $\Gamma$- {\it group}, a group $G$ together with another group $\Gamma$, an {\it operator group}, that acts on $G$ by automorphisms. Such an action is described by a homomorphism $\Gamma \to {\rm Aut}(G)$. We denote the action by $a^\g$, $\g \in \Gamma$ and $a \in G$.  An extended Alexander group is a  ${\Bbb Z}^2$-group. 

A {\it homomorphism} from a $\Gamma$-group $G$ to a $\Gamma'$-group $G'$ is a homomorphism $f:G \to G'$ of the underlying groups together with a homomorphism $\phi: \Gamma \to \Gamma'$ of operator groups  that intertwines the actions; that is,  $f(a^\g) = f(a)^{\phi(\g)}$ for every $\g \in \Gamma$ and $a \in G$. It is a {\it $\Gamma$-group isomorphism} if both $f$ and $\phi$ are group isomorphisms.  

Oriented virtual knot diagrams that are equivalent under generalized Reidemeister moves determine extended Alexander groups that are isomorphic in this sense. In fact, a ${\Bbb Z}^2$-group isomorphism can be found such that $\phi$ is the identity map.

The adjective ``extended" is a reference to the variable $v$. Setting $v$ equal to $1$ produces a ${\Bbb Z}$-group ${\cal A}_k$ called the {\it Alexander group} [SW01].
Regarded as an ordinary group, ${\cal A}_k$ is closely related to the commutator subgroup of $\pi_k$.\bs

\ni {\bf Proposition 2.1.} [SW01] The Alexander group ${\cal A}_k$ of an oriented virtual knot $k$ is isomorphic (as an ordinary group) to $\pi_k'*F$, the free product of the commutator subgroup of the group $\pi_k$ of the knot and a free group $F$ of countable rank.  Moreover, $F$ has a natural ${\Bbb Z}$-action, and ${\cal A} \cong \pi_k'*F$ as ${\Bbb Z}$-groups.\bs

Assume that $k$ is a classical knot. Let $p:\tilde X \to X$ be the universal abelian cover of the exterior $X$ of $k$. The term Alexander group is justified by the following.\bs

\ni {\bf Corollary 2.2.} [SW01] If $k$ is a classical knot, then the abelianization of ${\cal A}_k$, regarded as a ${\Bbb Z}$-module, is the Alexander module $H_1(\tilde X, p^{-1}(*); {\Bbb Z})$. \bs

A long virtual knot diagram acquires an orientation from ${\Bbb R}$. An extended Alexander group is defined as for virtual knots, and similar arguments show that equivalent long virtual knot diagrams produce isomorphic groups. 
\bs

\ni {\bf Remark 2.3.} Extended Alexander groups are defined in [SW01] more generally for any link $\ell$. The variable $u$ is replaced by $u_1, \ldots, u_d$ corresponding to the components of $\ell$, and general versions of the first part of Proposition 2.1 and Corollary 2.2 are valid.  Long virtual links can be defined, and  multivariable versions of extended Alexander groups can be defined for them. However, we confine our attention here to the case of knots. \bs

Given a diagram $D$ for an oriented long virtual knot $k$, let $a_{-\infty}$ (resp. $a_{+\infty}$) denote the generator family corresponding to the arc that runs to $-\infty$ (resp. $+\infty$). \bs

A $\Gamma$-{\it group system} is a tuple $(G; a_1, \ldots, a_r)$ consisting of a $\Gamma$-group $G$ and elements $a_i \in G$. Two $\Gamma$-group systems $(G; a_1, \ldots, a_r)$ and $(G'; a_1, \ldots, a'_r)$ are {\it isomorphic} if there is a $\Gamma$-group isomorphism $f: G \to G'$ such that $\phi: \Gamma \to \Gamma$ is the identity and such that $f(a_i) = a_i',\ $ for each $i$. It is immediate that the following ${\Bbb Z}^2$-group system is an invariant of long virtual knots. \bs

\ni {\bf Definition 2.4.} Let $k$ be a long virtual knot. The {\it extended Alexander group system of $k$}  is the triple $(\tilde {\cal A}_k; a_{-\infty}, a_{+\infty})$.  \bs

When $k$ is classical, the elements $a_{-\infty}$ and $a_{+\infty}$ are equal. This follows from the well-known fact that for a classical knot (or link) diagram, any one Wirtinger relation is a consequence of the the others. This does not apply to virtual knots, and $a_{-\infty}$ and $a_{+\infty}$ can be different. This will be seen in the examples below.  \bs


\ni {\bf 3. Examples.} (i) Consider the three long virtual knots $k_1, k_2$ and $k_3$ shown with extended Alexander group generators in Figure 2. If we join ends, the virtual knots obtained are the same, the {\it virtual trefoil} $k$.

The extended Alexander group $\tilde {\cal A}_{k_1}$ has defining relations 
$ac^u = db^u, a^v=b, db^u = ce^u, d^v =e.$ Using the second, third and fourth relations, we eliminate $b, c$ and $e$. The extended Alexander group of $k_1$ has presentation  $\langle a, d \mid ad^ua^{u^2 v}=d a^{uv}d^{u^2v}\rangle$. The extended Alexander group system is $(\tilde {\cal A}_{k_1}; a, d^v)$, which is isomorphic to  $(\tilde {\cal A}_{k_1}; a, d)$.

By similar calculations we find that the extended Alexander group systems of $k_2$ and $k_3$ are $(\tilde {\cal A}_{k_2}; a, d)$ and 
$(\tilde {\cal A}_{k_3}; a, dd^{uv^2} \bar d^{uv})$, where $\bar{}$ denotes inversion, and 

$$\tilde {\cal A}_{k_2}\cong \langle a, d \mid dd^{uv^2}d^{{u^2}v}=d^vd^ua^{{u^2}v}\rangle,\quad  \tilde {\cal A}_{k_3}\cong \langle a, d \mid dd^{uv^2}=d^va^{uv^2} \rangle.$$

An isomorphism between two of the group systems would map $a$ to itself. We see that no two group systems are isomorphic by considering the quotient groups $\tilde {\cal A}_{k_i}/\langle\langle a \rangle\rangle$, $i=1,2,3$, with respective presentations:
$$\langle d \mid dd^{u^2v} = d^u\rangle, \quad \langle d \mid dd^{uv^2}d^{u^2v}=d^vd^u \rangle, \quad \langle d\mid dd^{uv^2}=d^v\rangle.$$
The $0$th characteristic polynomials of the abelianized groups, regarded as modules over ${\Bbb Z}[u^{\pm 1}, v^{\pm 1}]$, are $u^2 v - u +1,\ u^2v + uv^2 -u - v +1$ and $uv^2 -v +1$, respectively. Since the polynomials, which are well defined up to multiplication by monomials $\pm u^iv^j$, are group invariants, $k_1, k_2$ and $k_3$ are distinct.\bs

\epsfxsize=3.5truein
\centerline{\epsfbox{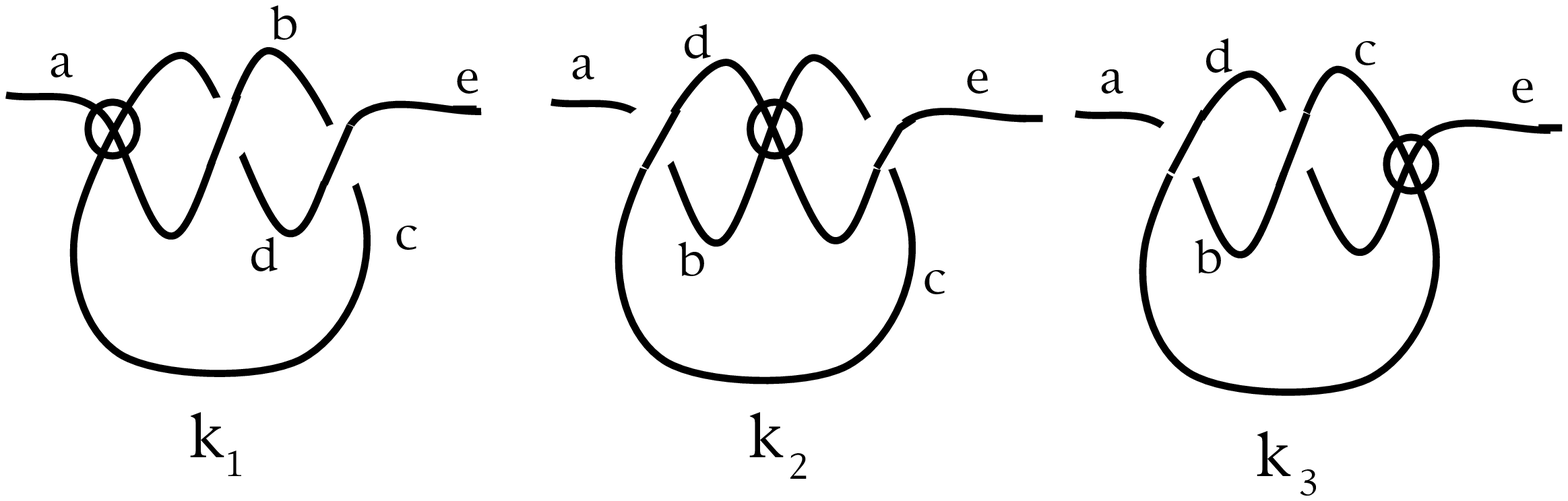}} \bs
\centerline{{\bf Figure 2:} Long virtual knots $k_1, k_2$ and $k_3$}\bs

(ii) Consider now the long virtual knot $k_4$ shown in Figure 2 with generators $a,b,c,d,e$ indicated. Its extended Alexander group $\tilde {\cal A}_{k_4}$ has relations
$a c^u = db^u,  d^v = c, bd^u=ec^u, b^v = c.$ As in the previous example, we use the second and fourth relations to eliminate $b$ and $c$, obtaining:

$$\tilde{\cal A}_{k_4} \cong \langle a,d,e\mid a=dd^u\bar d^{uv} =e\rangle$$

Obtain $k_5$ from $k_4$ by reversing its two classical crossings (Figure 3). Its extended Alexander group has relations $xz^u=ry^u, x^v=y, yr^u=sz^u, s^v =r.$ As in the previous example, we use the second and fourth relations to eliminate generators. We obtain: 
$$\tilde {\cal A}_{k_5} \cong \langle x,z,s\mid xz^u=s^vx^{uv}, x^vs^{uv}=sz^u\rangle.$$
By considering the quotients, we see that $k_4$ and $k_5$ and are distinct. The reader can check that $\tilde {\cal A}_{k_4}/\langle\langle a,e\rangle\rangle$ is nontrivial while 
$\tilde {\cal A}_{k_5}/\langle\langle x,s\rangle\rangle$ is trivial. \bs

\epsfxsize=2.5truein
\centerline{\epsfbox{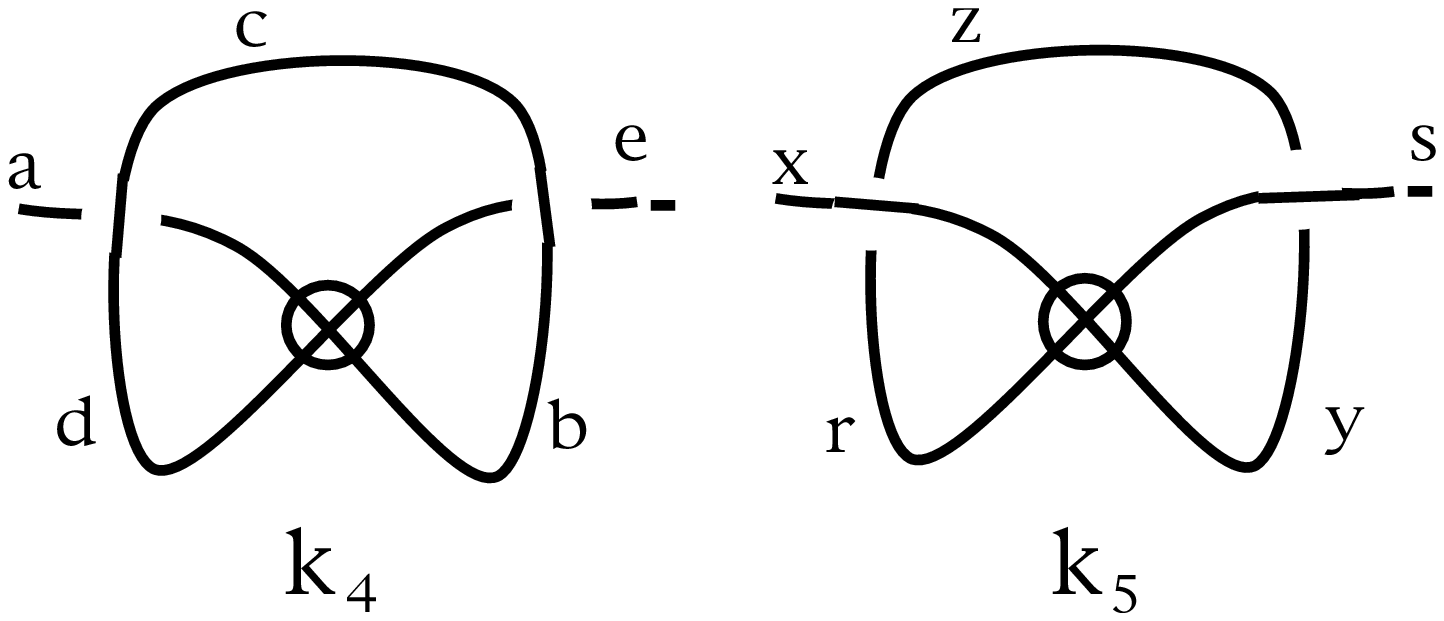}} \bs
\centerline{{\bf Figure 2:} Long virtual knots $k_4$ and $k_5$}\bs

(iii) Consider the two product knots $k_4\cdot k_5$ and $k_5\cdot k_4$. 
These were shown to be distinct by Manturov [M04] by counting homomorphisms of the associated quandle modules into ${\Bbb Z}_{43681}$. We establish the result by Alexander groups.

Using the notation for generators already established, it is a straightforward matter to check that the extended Alexander group of $k_4\cdot k_5$ has the following presentation.
$$\tilde {\cal A}_{k_4\cdot k_5} \cong \langle a,d,x,z,s \mid a = dd^u\bar d^{uv}=x, xz^u = s^vx^{uv}, x^vs^{uv}=sz^u\rangle.$$ 
Consequently, 
$$\tilde {\cal A}_{k_4\cdot k_5}/\langle\langle a\rangle\rangle \cong \langle d,z,s\mid dd^u=d^{uv}, z^u = s^v, s^{uv} = s z^u\rangle$$
$$\cong \langle d,s\mid dd^u = d^{uv}, s^{uv}=ss^v\rangle.$$
Similar computations show:
$$\tilde {\cal A}_{k_5\cdot k_4} \cong \langle x,z,a,d,e \mid xz^u = a^vx^{uv}, x^va^{uv}=az^u, a = dd^u\bar  d^{uv}=e\rangle,$$
and 
$$\tilde {\cal A}_{k_5\cdot k_4}/\langle\langle x\rangle\rangle \cong \langle z,a,d,e\mid z^u=a^v, a^{uv}=az^u, a=dd^u\bar d^{uv}=e\rangle.$$
We use the first relation to eliminate $z$, obtaining
$$\tilde {\cal A}_{k_5\cdot k_4}/\langle\langle x\rangle\rangle \cong \langle a,d\mid a^{uv}=a a^v, a = dd^u\bar d^{uv}\rangle$$
$$\cong \langle d \mid (dd^u \bar d^{uv})^{uv}=(dd^u\bar d^{uv})(dd^u\bar d^{uv})^v\rangle.$$
To see that $\tilde {\cal A}_{k_4\cdot k_5}/\langle\langle a\rangle\rangle$ and $\tilde {\cal A}_{k_5\cdot k_4}/\langle\langle x\rangle\rangle$ are not isomorphic, consider their abelianizations regarded as modules over 
${\Bbb Z}[t^{\pm 1}]$, where $u=v=t$. In the first case, the module is a direct sum of two cyclic modules, each with $0$th characteristic polynomial equal to $t^2-t-1$. In the second case, the module is cyclic.\bs

\ni {\bf 4. Colorings.} The well-known method of $p$-coloring knot diagrams, introduced by R.H. Fox in the early 1960's, extends to a method for distinguishing many long virtual knots. We explain its relationship to the Alexander group, and apply it to show that any virtual knot arises from infinitely many long virtual knots by making ends meet. 

Throughout, $p$ will denote an integer greater than $1$. Let $D$ be a classical knot diagram. A {\it $p$-coloring} of $D$ is an assignment of ``colors"  $\a, \b, \g, \ldots\in {\Bbb Z}/p{\Bbb Z}$ to the arcs of $D$ so that at every crossing twice the value assigned to the overcrossing is equal to the sum of the values assigned to the two undercrossing arcs.  Every diagram admits $p$ ``monochromatic" colorings. By a {\it nontrivial $p$-coloring} we will mean a $p$-coloring using more than one color. 

It is an elementary fact that different diagrams for $k$ admit the same number of $p$-colorings. Much more is known.  The set of $p$-colorings for a diagram form a ${\Bbb Z}/p{\Bbb Z}$-module in a natural way, and this vector space is isomorphic to $H^1(M_2; {\Bbb Z}/p{\Bbb Z}) \oplus {\Bbb Z}/p{\Bbb Z}$, where $M_2$ denotes the $2$-fold cyclic cover of $S^3$ branched over $k$. However, we will not require this fact. 

Fox's interest in $p$-colorings came from considerations of the knot group $\pi_k$.  The set of $p$-colorings of a diagram for $k$ are in one-to-one correspondence with the set of homomorphisms from $\pi_k$ onto the dihedral group $D_p$ of order $2p$. The idea of $p$-colorings can be found in Chapter 10 of [CF63].

Alternatively, one can regard the commutator subgroup $\pi_k'$ as a ${\Bbb Z}$-group, the action of a generator ${\Bbb Z}$ being the automorphism of $\pi_k'$ induced by conjugation by a meridian. The cyclic group ${\Bbb Z}/p{\Bbb Z}$ is a ${\Bbb Z}/2{\Bbb Z}$-group, with the action of the generator of ${\Bbb Z}/2{\Bbb Z}$ taking any element of  ${\Bbb Z}/p{\Bbb Z}$ to its inverse. (The dihedral group $D_p$ is the semidirect product of ${\Bbb Z}/p{\Bbb Z}$ by ${\Bbb Z}/2{\Bbb Z}$ with this action.) There is a bijective correspondence between nontrivial $p$-colorings of a diagram for $k$  
and nonzero homomorphisms from $\pi_k'$ to ${\Bbb Z}/p{\Bbb Z}$ with operator groups ${\Bbb Z}$ and ${\Bbb Z}/2{\Bbb Z}$ mapping via natural projection.\bs
 
 \ni {\bf Remark 4.1.} See [SW98] for a different interpretation of Fox $p$-colorings. In this setting $p$-colorings correspond to period-2 points of a symbolic dynamical system associated to $k$. In general points of period $r$ correspond to elements 
 of $H^1(M_r; {\Bbb Z}/p{\Bbb Z}$), where $M_r$ is the $r$-fold cyclic cover of $S^3$ branched over $k$. \bs  
 
Fox $p$-colorings have been extended previously for virtual knots (see for example [KSW00], [M04]). The theory extends further to long virtual knots.
The set of $p$-colorings of a diagram has the structure of a ${\Bbb Z}/p{\Bbb Z}$-module. One can examine the effect of generalized Reidemeister moves to see that equivalent diagrams determine isomorphic modules.  
Alternatively, one can appeal to the following, the proof of which is straightforward. \bs

\ni {\bf Proposition 4.2.} Let $D$ be a diagram of a long virtual knot $k$. The module of $p$-colorings of $D$ is isomorphic to ${\rm Hom}({\cal A}_k, {\Bbb Z}/p{\Bbb Z})$, where ${\cal A}_k$ (resp. ${\Bbb Z}/p{\Bbb Z}$) is regarded as a ${\Bbb Z}$-group (resp. ${\Bbb Z}/2{\Bbb Z}$-group) and the map between operator groups is mod-$2$ reduction. \bs

The abelianization of the Alexander group ${\cal A}_k$ of a long virtual knot is a ${\Bbb Z}$-module. It can be described by
${\Bbb Z}^N/{\Bbb Z}^{N-1} A(t)$,
where $N$ is the number of arcs of a diagram for $k$ and $A(t)$ is an $(N-1)\times N$ presentation matrix with entries in ${\Bbb Z}[t^{\pm 1}]$. \bs

\ni {\bf Definition 4.3.} The {\it determinant} of $k$ is the greatest common divisor of all $(N-1)\times (N-1)$ minors of $A(-1)$. \bs

Definition 4.3 is consistent with the definition of determinant for classical knots regarded as closings of long knots (see \S1). It is well known that a diagram for a classical knot has a nontrivial $p$-coloring, for $p$ prime, if and only if $p$ divides the determinant of $k$. The proof requires only standard linear algebra. The same argument shows: \bs

\ni {\bf Proposition 4.4.} A long virtual knot $k$ has a nontrivial $p$-coloring if $p$ divides the determinant of $k$. The converse is true whenever $p$ is prime. \bs

The following result is well known for classical knots.\bs

\ni {\bf Proposition 4.5.} The determinant of any long virtual knot is odd. In particular, it is nonzero. \bs

\ni {\bf Proof.} Consider a long virtual knot $k$ with diagram $D$ as above. The determinant of $k$ is odd if and only if some $(N-1)\times (N-1)$ minor of $A(-1)$ is odd. Since the entries of $A(-1)$ are congruent modulo 2 to those of $A(1)$, it suffices to prove that some $(N-1)\times (N-1)$ minor of $A(1)$ is odd. 

The matrix $A(1)$ presents the ${\Bbb Z}$-module with generators $a,b,c,\ldots$ corresponding to the arcs of $D$ and relations corresponding to each crossing, identifying each of the two undercrossing arcs. Since $k$ is a knot, all of the generators are identified, and hence the module is isomorphic to ${\Bbb Z}$. Consequently, each $(N-1)\times (N-1)$ minor of $A(1)$ is a unit $\pm 1$. \qed
\bs

\ni {\bf Example 4.6.} Let $D$ be a diagram for a long virtual knot that closes to a virtual knot $k$. Construct a sequence $D_n$ of diagrams for long knots $k_n,\ n \in {\Bbb N}$, as indicated in Figure 4. $D_{n+1}$ is obtained from $D_n$ by adding a classical crossing followed by a virtual crossing above and below $D$ in such a way that the classical crossings alternate in sign. Note that $D$ is embedded in the diagrams $D_n$. The reader can verify that by joining the ends of any $D_n$, we obtain a diagram for $k$. Hence each long virtual knot $k_n$ maps to $k$ under the closing map $\kappa$ (see \S1). 

\epsfxsize=1.5truein
\centerline{\epsfbox{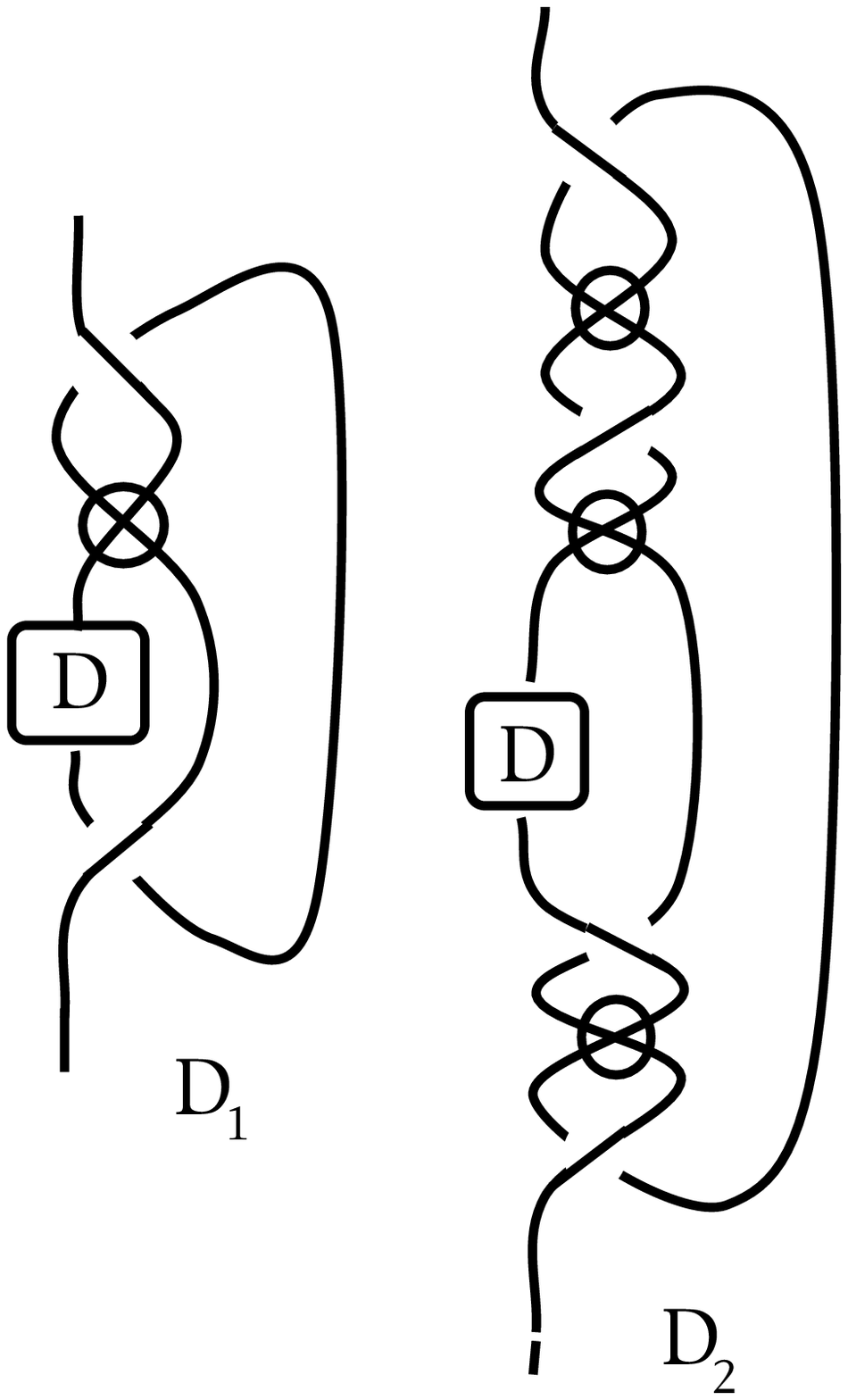}} \bs
\centerline{{\bf Figure 4:} Long virtual knot diagrams $D_n$}\bs

We will use $p$-colorings to show $\{k_n\}$ contains infinitely many distinct long virtual knots. Color the top two arcs of any diagram $D_n$ with colors $\a$ and $\b$ in ${\Bbb Z}/p{\Bbb Z}$. The colors of the arcs above the level of $D$ are uniquely determined. Let $\g$ be the color entering $D$. {\sl Color every arc of $D$ with color $\g$.} Now the colors of arcs below the level of $D$ are also determined. A $p$-coloring of $D_n$ is obtained if and only if the color of the arc emerging at the bottom right-hand side of the diagram is equal to $\b$. We leave it to the reader to check that nontrivial $p$-colorings exist if and only if there extist
distinct $\a$ and $\b$ such that 

$$(\a, \b)\cdot S(TS)^{n-1}U = (\delta, \b),\ {\rm for\ some}\ \delta$$
where 
$$S= \pmatrix{1& 2\cr 0&-1}, \quad T= \pmatrix{-1& 0\cr 2&1}, \quad U= \pmatrix{0& -1\cr 1&2}.$$
By a straightforward induction argument, the condition amounts to 
$$(2n+1)(\a- \b) =0.$$
Hence $D_n$ admits a nontrivial $(2n+1)$-coloring. By Proposition 4.4, if $2n+1$ is prime, then it divides the determinant of $D_n$. In view of Proposition 4.5, the collection $\{k_n\}$ is infinite. \bs

\ni {\bf References.} \bs

\ni [CF63] R.H. Crowell and R.H. Fox, Introduction to Knot Theory, Ginn and Co., Boston, 1963 \ss

\ni [GPV00] M. Goussarov, M. Polyak and O. Viro,  Finite type invariants of classical and virtual knots, {\sl Topology\ \bf 39} (2000), 1045--1068. \ss

\ni [K97] L.H. Kauffman, Talks at MSRI Meeting (January 1997), AMS Meeting at University of Maryland (March 1997), Isaac Newton Institute Lecture (November 1997), Knots in Hellas Meeting, Delphi, Greece (July 1998), APCTP-NANKAI Symposium on Yang-Baxter Systems, Non-linear Models and Applications, Seoul, Korea (October 1998). \ss

\ni [K99] L.H. Kauffman, Virtual knot theory, {\sl European J. Comb.\ \bf 20} (1999), 663--690. \ss

\ni [KSW00] D.A. Krebes, D.S. Silver and S.G. Williams, Persistent invariants of tangles, {\sl J.\ Knot\ Theory\ and\ its\ Ramifications\ \bf 9} (2000), 471--477.\ss

\ni [M04] V. Manturov, Long virtual knots and their invariants II: the commutation problem, preprint. \ss

\ni [SW98] D.S. Silver and S.G. Williams, Generalized $n$-colorings of links, 
{\sl Knot Theory}, Banach Center Publ. {\bf 42}. Warsaw 1998. \ss

\ni [SW01] D.S. Silver and S.G. Williams, Alexander groups and virtual links, {\sl J. Knot Theory and its Ramifications\ \bf 10} (2001), 151--160. \ss

\ni [SW03] D.S. Silver and S.G. Williams, Polynomial invariants of virtual links, 
{\sl J. Knot Theory and its Ramifications\ \bf 12} (2003), 987--1000. \ss

\end